\documentclass[11pt]{aaa-art}

\advance\textheight0.4cm
\usepackage{amssymb}
\usepackage{mathrsfs}
\renewcommand{\P}{{\mathscr P}}
\newcommand{\R}{{\mathbb R}}
\newcommand{\N}{{\mathbb N}}

\newcommand{\opc}{o.p.c.}

\newcommand {\x}{{\tt x}}

\renewcommand{\O}{{ O}}

\def \itm#1 {\item[(#1)]}

\swapnumbers
\theoremstyle{remark}

\def\xx#1 {\newtheorem{#1}[thm]{#1}}
\xx Lemma
\xx Theorem
\xx Corollary
\xx Idea
\xx {First Step}
\xx Consequence
\xx Example
\xx Definition
\xx Problem
\xx{Open Questions}
\xx{Open Question}
\xx Speculation

\author{Martin Goldstern}
\title{Lattices, Interpolation, and Set Theory}
\address{Institut f\"ur Algebra\\TU Wien\\
Wiedner Hauptstr 8-10/118.2\\A-1040 Wien}
\email{Martin.Goldstern@tuwien.ac.at}
\urladdr{http://info.tuwien.ac.at/goldstern/}

\thanks{This paper is available on my homepage, and also from {\tt
arXiv.org}}

\thanks{I am grateful to the Austrian Science Foundation (FWF) for
supporting my research through grant P13325-MAT}

\begin{document}
\nocite{mfl}
\nocite{l01}

\nocite{633}
\nocite{Haviar+Ploscica:1998}
\nocite{Wille:1977}

\begin{abstract}
We review a few results concerning interpolation of monotone functions
on infinite lattices, emphasizing the role of set-theoretic
considerations. 
We also discuss a few open problems. 
\end{abstract}

\maketitle
\section{Introduction}
\label{section.intro}

In this article I will present  a few ideas concerning the
relationship between lattice polynomials and monotone functions 
on lattices.    It turns out that (for infinite lattices) there are 
purely algebraic questions which (apparently) cannot be solved with 
purely algebraic methods;  ideas from logic --- in
particular, from set theory  --- have to be
used.  Most prominent among these logical concepts is the notion of
distinct infinite cardinalities, but also the more subtle notion of
``cofinality'' plays a role. 

Most of the results here  are not new; specifically, sections 
\ref{section.opc}
 and \ref{section.IP}
 just
highlight  ideas from papers that are already published. 
  The reason for including those known
results here is that they serve well to show  the interconnections
between algebraic and set theoretic ideas. \\
 The observations in sections 
\ref{section.complete},  \ref{section.sigma}, 
and \ref{section.ortho} are based on these 
older results,  but are themselves new.

I will often not give complete proofs but treat only a characteristic
case, in the hope that the information presented here is enough for
the reader to complete the proof. 
For example, I 
  will sometimes prove (or even state) a theorem which is true for $n$-ary
 functions only for unary functions.     The general case is often just
an easy generalization (but sometimes needs an additional idea). 
 In  any case, complete proofs are
available elsewhere in the literature.

To avoid some case distinctions, we
 will only consider bounded lattices; in addition to the operations
 $\vee $ and~$\wedge$, all lattices that we consider in this paper 
come with constants $0$ and~$1$.  So ``lattice'' will mean ``bounded
 lattice'', all homomorphisms 
 have to respect 0 and 1, all sublattices contain $0$ and~$1$, etc.

\section{Order-polynomial completeness}
\label{section.opc}

The problem that got me interested in the subject was the following:
\begin{quote}
It is clear that every lattice polynomial\footnote{Recall that the set
$L[\x_1,\ldots, \x_n]$ of $n$-ary lattice polynomials [``with
coefficients from $L$''] is the free product of $L$ and the free lattice 
with $n$ generators, or in other words: take 
 the set of well-formed expression
using the indeterminates $\x_1,\ldots,\x_n$, elements of the
lattice~$L$ and the operation symbols $\vee $ and~$\wedge$, and 
divide by
the obvious equivalence relation: ``equal in all extensions of $L$''.
 The $n$-ary 
 function induced by a lattice polynomial is called ``polynomial
function'' (or ``algebraic function'' by some authors).}  
 represents a monotone
function.  \\For which lattices is the converse true? 
\end{quote}

Let us call a lattice $L$ $k$-\opc\ ($k$-order polynomially complete) if
every monotone function $f:L^ k \to L$ is represented by a polynomial,
and call $L$ \opc\ if $L$ is $k$-\opc\ for every~$k$. 

The finite \opc\  lattices were characterized by Wille (see section
\ref{section.finite}). 
In this section  we will show the main ideas  of the proof 
of the following theorem: 

\begin{Theorem}
\label{opc}
There are no infinite \opc\  lattices. 
\end{Theorem}

Although the problem itself 
is now solved, I feel that our understanding is
not at all complete.   As evidence for this lack of
understanding, see problems~\ref{A}  and \ref{B}. 

\subsection*{Cardinalities}    If there are ``more'' monotone
functions than polynomials, then clearly not all monotone
functions are polynomials.  In the finite case this does not help much,
since there are only finitely many monotone functions for every fixed
arity,  while there are always  infinitely many 3-ary polynomials. 
{\smaller\relax [Of course, for lattices with extra properties, such as
distributive lattices, we can get a normal form for polynomials, which
means that 
the set of polynomials looks  rather ``small'' to us.]}

In infinite lattices the situation is reversed:  The set $P$ of
polynomials on a lattice $L$ has the same
 cardinality\footnote{Recall 
the following elementary facts about cardinals:
\begin{itemize}
\item 
 We write $A \approx B$ or $|A|=|B|$ iff there is a
bijection between $A$ and~$B$.  We write $ |A|$ for the
``cardinality'' of~$A$, which can informally be defined as the
equivalence class of $A$ with respect to~$\approx$. 
\item  We write $A \le  B$ or $|A|\le|B|$ iff there is a 1-1 map from 
 $A$ into~$B$.   Equivalently (for nonempty $A$):  $|A|\le|B|$
 iff $B$ can be mapped onto~$A$. 
\item The notation   $|A|\le|B|$ is justified by the following
theorem: If   $|A|\le|B|$ and  $|B|\le|A|$, then $|A| = |B|$. 
\item $|A| = |A^ n | = |\bigcup_k A^ k|$ for all infinite~$A$. 
\item We say that a set $A$ is countable 
 and write $|A| = {\aleph_0} $
iff there is a bijective map from $\N$ onto~$A$. $A$ is ``at most
 countable''  [$|A|\le {\aleph_0} $] iff $A$ is countable or finite. 
\end{itemize}
}
 as $L$ itself,
whereas the set $L^ L$ of all functions from $L$ into $L$ has always
strictly larger cardinality.

However, not all functions are monotone, and there are even trivial
examples of lattices $L$ with only ``few'' monotone functions, i.e.,
satisfying 
$$ \{f:  \mbox{$f$ is monotone from $L$ to $L$}\}  \approx L$$
for example the linear order $\R$ (the real numbers).

Here is the main idea that can be used for many lattices $L$ to show
that $L$ is not \opc: 

\begin{enumerate}
\itm a find a ``nice'' substructure $A$ of the same cardinality as $L$
\itm b find a subset  $L' \subseteq L$
  such that 
\begin{itemize}
 \item  for  all ``nice'' structures $A$ of 
cardinality $ \kappa$  there are  $> \kappa $ many monotone maps from
            $A$ to $L'$  
\item the set $L'$, with the order inherited from $L$, is a complete
lattice [not necessarily a sublattice of $L$]
\end{itemize}
\itm  c conclude that there are  $> \kappa $ many monotone maps from
$L$ to $L$.  
\end{enumerate}

Part (c) is easy: every monotone $f:A \to L'$ can be extended to $\bar
f:L \to L'$ via 
$$ \bar f(x) = \sup \nolimits_{L'}\{ f(z): z\in A, z\le x\}.$$


\subsection*{From (anti)chains to monotone functions}

It is well-known (see, e.g., \cite{Haviar+Ploscica:1998} or 
\cite{633}) that  if $A \subseteq L $ is an antichain\footnote{$A
\subseteq L$ is 
called an antichain if no two distinct elements of $A$ are
comparable.} then there are\footnote {If $|A| = \kappa$, we
write $2^ \kappa $ for the cardinality of~$\P(A)$.  Thus, a set $B$ is
of cardinality $2 ^ {\aleph_0} $ iff $B$ can be mapped bijectively
onto $\P(\N)$ or~$\R$.}
 $2^ {|A|} > |A|$ many pairwise incomparable\footnote{If 
 $f,g: L^ n\to L$ are
functions, we say $f\le g$ iff
 $f(a_1,\ldots, a_n) \le g(a_1,\ldots, a_n)$
 holds for all $(a_1, \ldots, a_n) \in L^ n$.
 The relation $\le$ is then  a partial order of functions.    
The polynomials are (via
the functions they induce) quasiordered. 
}
 functions from $A$ to~$\{0,1\}$, and since $\{0,1\}$ is a complete
lattice  it is trivial to extend all those functions to 
monotone functions defined on all of $L$. If 
 $A \subseteq L$ is
well-ordered\footnote{Recall that 
a linearly ordered set $(A,\le)$ is called ``well-ordered'' iff every
nonempty subset of $A$ has a least element.}, there are $2^{|A|}$
many monotone functions from $A$ to~$A$, and again they can all be
extended to total monotone functions from $L$ 
into the complete lattice $A \cup
\{0,1\}$.  

Thus, a first approximation to the program (a)(b)(c) outlined above
is:    For a given lattice~$L$, find a ``large'' antichain or
well-ordered chain in $L$ (where ``large'' means: of the same
cardinality as $L$ itself). 

Unfortunately, this is not always possible but at least for countable
lattices the following observation helps: 
{\em Ramsey's Theorem}\footnote{We will use the following 
version of the infinitary Ramsey theorem:
  Whenever the edges of the complete graph on countably
many vertices are colored with 3 colors, then there is an infinite
complete subgraph, all of whose edges have the same color.}
 implies that every infinite partial order
contains either a chain isomorphic or anti-isomorphic to the natural
numbers, or an infinite antichain. 	

{\small Proof:  We may assume that the partial order $P$ is countable:
$P = \{p_1, p_2, \ldots \}$, where all $p_i$ are distinct.   Color the
edges of the complete graph on $\{1,2,\ldots\}$ with 3 colors:  The
edge $\{i,j\}$ is colored according to whether the map
$i \mapsto p_i$,
$j\mapsto p_j$ is an isomorphism, anti-isomorphism, or neither.   Any
infinite complete
 subgraph whose edges are  colored with only one color will be
the desired chain or antichain.}

Hence any infinite partial order, in particular any infinite lattice,
admits $\ge 2^ {\aleph_0} $ many
 monotone functions, so an \opc\  lattice has to have size 
$\ge 2^{\aleph_0} $.

\bigskip

The above method yielding many incomparable  functions from an
antichain is set-theoretical.  The next idea, converting a set of
incomparable functions into an antichain, is  mainly algebraic. Here 
 it is important that we are interested in fully \opc
(and not only $1$-\opc) lattices: 

\subsection*{From monotone functions to antichains}
\begin{Lemma}\label{hp}
 Let $L$ be a lattice on which there are $\kappa $ many
pairwise  incomparable
 $k$-ary polynomials, 
where\footnote {$cf(\kappa) > {\aleph_0} $ 
[read: $\kappa $ has uncountable cofinality] 
 means:  $\kappa $ is the
cardinality of an infinite  set $A$ with the following property: 
\begin{quote} Whenever $ A = \bigcup_{n=1}^ \infty 
 A_n$, then for some $n$, $|A| =
|A_n|$.
\end{quote}
  For example, it is true
 (but not quite trivial, unless we assume the continuum hypothesis)  
that the set $\R$ of real
numbers satisfies $cf(|\R|) > {\aleph_0} $.
\endgraf 
The negation of this property is denoted $cf ( \kappa) = {\aleph_0}$:
$\kappa = |A| $ for some infinite set $A$ which can be written 
 as a countable union of sets of strictly smaller cardinality.
}
 $cf(\kappa ) > {\aleph_0}$.

Then there is some $n'$ such that  $L^ {n'}$ contains an antichain of
size $ \kappa $. 
\end{Lemma}

\begin{proof} (This idea, and in fact the whole proof, is due to
\cite{Haviar+Ploscica:1998}.) 

Assume that the polynomials $(p_i(\bar \x): i \in I)$ are pairwise
incomparable, where $|I| = \kappa $, and $\bar \x$ abbreviates
$(\x_1,\ldots, \x_k)$.  
 For each $i$ there is a natural
number $n_i$ such that   $p_i( \bar \x)$ can be written as 
$t_i(\bar \x, \bar c_i)$, where $t_i$ is a $(k+n_i)$-ary term, and $\bar
c_i \in L^{n_i}$ (the entries of the vector $\bar c_i$ are the
``coefficients'' of the polynomial $p_i$). 
Since there  only countably many
terms, our assumption $cf( \kappa ) > {\aleph_0}$ implies we can thin out
the set $I$ to a set $I'$ of the same cardinality such that  all the
$(n_i,t_i)$ (for $i \in I'$) are equal, say $= (n', t')$.   Now it is
easy to check that the ``coefficients'' $(\bar c_i: i \in I')$ form an
antichain in~$L^{n'}$.  
\end{proof}

\begin{Corollary}\label{hp2}
 If $L^ {n_1}$ has an antichain or well-ordered chain 
of size~$\kappa$,
and $L$ is $n_1$-\opc, then there are $2^ \kappa $ many pairwise
incomparable polynomial functions $p:L^ {n_1} \to L$. So there is some
$n_2$  such that  $L^ {n_2} $ has an antichain of size $2^ \kappa$. 

Repeating this argument we can find some $n_3$ such that $L^ {n_3}$
has an antichain of size $2^{2^ \kappa }$, etc. 
\end{Corollary}

So let $L$ be infinite and  \opc{}  
To
  simplify some calculations, we will assume
 GCH, the generalized continuum hypothesis,  for the rest of this
 chapter.  However, all of what will be said will --- with some obvious
modifications --- remain true even without this additional
  assumption. 

We know that $L$ contains a chain or antichain of size~${\aleph_0}$,
so also one of size $2^ {\aleph_0} $ (since we assume GCH we can write
this as~${\aleph_1}$).   Iterating corollary~\ref{hp2} we get 
$$ |L| \ge {\aleph_1}, \ 
 |L| \ge {\aleph_2}, \  |L| \ge {\aleph_3}, \ldots$$
and finally $|L| \ge {\aleph_\omega } $.

It turns out that the proof splits into three very different
cases, depending on $\kappa$, the  cardinality of $L$:
\begin{description}
\item[{a}] For some $\mu < \kappa$, $2^\mu \ge \kappa$. 
\item[{b}] Case (a) does not hold, and $cf(\kappa) = {\aleph_0} $.
\item[{c}]  Case (a) does not hold, and $cf(\kappa) > {\aleph_0} $.
\end{description}

Case (a) is treated similarly to the cases $\kappa=\aleph_1$,
$\kappa=\aleph_2$, etc.  If case (a) does not hold, i.e., if we have 
$$ \forall \mu <\kappa: 2^\mu < \kappa$$
then we call $\kappa$ a ``strong limit cardinal''. 

We will in this paper repeat the argument from 
\cite{633} for case (b).  (Case (c) was considered in \cite{opc}.) 
It turns out that the case $|L| = {\aleph_\omega } $ is typical for
(b), so to again simplify the notation we will assume this equality. 

Thus $L$ can be written as $L =  L_0 \cup L_1 \cup \cdots\, $, where
$|L_n| = {\aleph_n} $,  $L_{n+1} \approx \P(L_n)$.

Can we prove that every lattice  $L$ of size ${\aleph_\omega}$
contains either an antichain or a chain (well-ordered or dually
well-ordered) of cardinality~${\aleph_\omega}$?   Unfortunately this
is not true.   The best we can do is to show that for each~$n$, $L$
must contain a (well-ordered or co-well-ordered) chain or antichain
$A_n$ of cardinality ${\aleph_n} $, but it is quite possible that the
union $A:= A_0 \cup A_1 \cup A_2 \cup \cdots$ is neither a chain nor
an antichain, as the following examples show: 

\begin{Example}
\label{ex1} For each $n$ let  $A_n$ be a well-ordered (or
co-well-ordered) set of cardinality ${\aleph_n} $, and assume that the
union 
$$L  = \{0\} \cup \{1\} \cup A_0 \cup A_1\cup  \cdots$$
is a disjoint union. Define a partial order on $L$ by requiring $0$
and $1$ to be the least and greatest elements respectively, and by
also requiring 
$$ \forall n \not= k \  \forall a \in A_n 
			\, \forall b \in A_k: 
\qquad \mbox{$a$ and $b$ are incomparable.}$$

Thus, $L$ consists of countably many chains (of increasing sizes),
arranged side-by-side.    We leave it to the reader to check that $L$
is indeed a lattice, in fact a complete lattice. 

$L$ contains chains of cardinality ${\aleph_n} $ for each~$n$,
but no chain of length~${\aleph_\omega}$, and every antichain in $L$
is at most  countable. 
\end{Example}

\begin{Example}
\label{ex2} For each $n$ let  $A_{2n+1}$ be an antichain of
cardinality~${\aleph_n}$, and let $A_{2n} = \{a_n\}$. 
Again, assume that the
union 
$$L  = \{0\} \cup \{1\} \cup A_0 \cup A_1\cup  \cdots$$
is a disjoint union. Define a partial order on $L$ by requiring $0$
and $1$ to be the least and greatest elements respectively, and by
also requiring 
$$ \forall n <  k \  \forall a \in A_n 
			\, \forall b \in A_k: 
\qquad a < b .$$

Thus, $L$ consists of countably many antichains, each one on top of
the previous one.   Again it is easy to check that $L$ is a complete
lattice.  [If $a,b\in A_{2n+1}$ are distinct, then $a\wedge b =
a_{2n}$.]

This time 
$L$ contains antichains of cardinality ${\aleph_n} $ for each~$n$,
but no antichain of length~${\aleph_\omega}$, and every chain in $L$
is at most  countable. 
\end{Example}

\begin{Example}
\label{ex3} For each $n$ let  $A_n$ be a well-ordered set of size $\aleph_n$. 
Again, assume that the
union 
$$L  = \{0\} \cup \{1\} \cup A_0 \cup A_1\cup  \cdots$$
is a disjoint union. Define a partial order on $L$ by requiring $0$
and $1$ to be the least and greatest elements respectively, and by
also requiring 
$$ \forall n <  k \  \forall a \in A_n 
			\, \forall b \in A_k: 
\qquad a > b .$$

Thus, elements of $A_0$ are above all elements of $A_1$, etc. $L$
itself is a chain of cardinality $\aleph_\omega$, but $L$ is not
well-ordered.  Moreover, any well-ordered subset of $L$ has
cardinality $< \aleph_\omega$, and any co-well-ordered subset of $L$ is
countable.  There are no antichains of size $>1$ in $L$.
\end{Example}

In  the first two 
examples we have a lattice $L$ of cardinality~${\aleph_\omega}
$, all of whose chains and antichains are of size~$<
{\aleph_\omega}$.     

Still, it is easy to see that the lattices in all three examples 
admit  $2^ {\aleph_\omega} $ many monotone
unary functions.  E.g., in example~\ref{ex2}, any map
$f:L \to L$ satisfying $f[A_n] \subseteq A_n$ for all~$n$, $f(0)=0$,
$f(1)= 1$  will be monotone, and there are
 $2^ {\aleph_0}  \times 2^ {\aleph_1}  \times \cdots \ = 2^
{\aleph_\omega }  $ many such maps. 

Hence, the following theorem suffices to prove that a lattice $L$ of
cardinality ${\aleph_\omega} $ cannot be \opc: 

\begin{Theorem}
If $L$ is a partial order of cardinality~${\aleph_\omega}$, then $L$
either contains a sufficiently large antichain or 
 (well-ordered or co-well-ordered) chain, or 
 a partial order $P$ which is isomorphic or antiisomorphic to one of
those  constructed in examples \ref{ex1},  \ref{ex2} and \ref{ex3}.
\end{Theorem}
This theorem can be easily deduced from the ``canonization'' theorem
of Erd\H os, Hajnal and Rado, see
\cite[28.1]{Erdos+Hajnal+Mate+Rado:1984}.

\goodbreak

The first ``open'' question is rather ill-defined: 
\begin{Problem}\label{A} Find a purely  algebraic proof that
every  \opc\  lattice must be finite. 
\end{Problem}

Why do we need such a proof?  First, there may be algebraists who are
 uncomfortable with the notions ``cardinality'' and ``cofinality''
 which were used in the proof above. 
 Secondly, and more to the point, a new proof of theorem~\ref{opc} may
 also shed light on the following problem, which is still 
open: 

\begin{Problem}\label{B}
Can there be an infinite   1-\opc\  lattice? If yes,
what about 2-\opc? etc. 
\end{Problem}

\begin{Speculation}
In \cite{opc} it is shown  that theorem \ref{opc} cannot be proved
without some weak version of AC, the axiom of choice; this seems to
indicate that some set-theoretic sophistication is necessary for any
proof of theorem \ref{opc}.  However, on closer scrutiny it turns out
that the only version of AC that was shown to be necessary for such a
proof is a statement that is strictly weaker than 
\begin{quote}
``every infinite set contains a countable subset''
\end{quote}
which for most non-logicians is not even recognizable as a version of
AC, so a  ``purely algebraic'' proof might still be possible.  
\end{Speculation}

\section{Finite \opc\  lattices}
\label{section.finite}

For the investigation of finite lattices, set theory plays of course
no role. We give the following characterization only to contrast it
below  with the situation for infinite lattices. 

Let us call a function $L\to L$ ``regressive'' if 
$ \forall x\in L: f(x)\le x$. 

\begin{Definition}
We say that a lattice satisfies Wille's property   if the only
regressive join-homomorphisms are the identity map and the constant
0. 
\end{Definition}
\begin{Theorem}
A finite lattice is \opc\  iff it is simple (in the algebraic sense) 
 and satisfies Wille's property. 
\end{Theorem}
This theorem is proved in \cite{Wille:1977}.

 It is of course crucial for the proof of this theorem
that the lattice under consideration be  finite, since the length of
the constructed polynomial will typically increase with the size  of
the lattice. 

The following easy example shows that the characterization theorem
cannot work for infinite sets (of any cardinality).  

\begin{Example}
\label{example-a}
Let $A$ be any infinite set disjoint from~$\{0,1\}$.   Define a
lattice $M_A = A \cup \{0,1\}$ by requiring $A$ to be an antichain and
$0\le a\le 1$ for all
$a\in A$.   

 Then $M_A$ is simple and
satisfies Wille's property, but is of course not \opc.  
\end{Example}

\section{Interpolation} 
\label{section.IP}

There are several natural ways to generalize the question ``which
monotone functions are represented by polynomials'' from finite to
infinite lattices.   The first approach, the property \opc, 
 was discussed in section \ref{section.opc}.  A second approach is the
following:
\begin{Definition}
We say that a lattice  $L$ has the 
IP (interpolation property) iff
\begin{quote}
for every monotone function $f:L^ n\to L$ and every finite $A \subseteq
L^ n$ there is a lattice polynomial $p$ such that 
$$  \forall (a_1,\ldots, a_n)\in A: \ f(a_1,\ldots,a_n) =
p(a_1,\ldots,a_n).$$
\end{quote}
\end{Definition}

In other words, if we equip $L$ with the discrete topology, and
$L^{(L^n)}$, the set of all functions from $L^n$ to~$L$, with the
product topology (= Tychonoff topology = topology of pointwise
convergence), then the IP says: 
\begin{quote} For all~$n$, the set of monotone functions in
$L^{(L^n)}$ is the closure of the set of all polynomial functions.
\end{quote}

Clearly, a finite lattice has the IP iff it is \opc, but
this is not true for infinitary lattices.  For example, the lattice
given in example \ref{example-a} has the IP, but is of course not
\opc{}

While there are no infinite
\opc\  lattices, the following can be shown: 
\begin{Theorem}
For every  lattice $L$ there is an infinite lattice $\bar L$
with the IP, where  $L$ is a sublattice of~$\bar L$. 
\end{Theorem}
Elementary cardinal arithmetic shows that $\bar L$ can be chosen to be
of the same cardinality as the cardinality of~$L$, as long as $L$ is
infinite.  {\small (If $L$ is finite, then it is known that $L$ can be
extended to a finite \opc\  lattice.)}

This theorem can be proved with a combination of set-theoretic and
algebraic ideas.  The algebraic content of the theorem is
represented by the following lemma
 (which is really a
weak version of the theorem itself). 
\begin{Lemma}\label{extend}
If $L$ is a lattice, $f:L\to L$ a monotone partial function,
 then there is a lattice $L'$
extending~$L$, such that $f$ is the restriction of a polynomial
function with coefficients in~$L'$. 
Moreover, if $L$ is complete,
 then $L'$ can be chosen to be again a complete
lattice, which is moreover an ``end extension'' of~$L$, i.e., 
$$ \forall b \in L'\, \forall a \in L \setminus \{1\}:  \ 
\bigl(b \le a \Rightarrow b \in L\bigr).$$
\end{Lemma}
(This lemma is proved in \cite{mfl} and \cite{l01}.)

It turns out that this lemma is sufficient to construct lattices with
stronger interpolation properties.  For example:

\begin{Definition}
We say that a lattice  $L$ has the 
${\sigma}$-IP (countable interpolation property)  iff
\begin{quote}
For every monotone function $L^ n\to L$ and every at most countable 
$A \subseteq
L^ n$ there is a lattice polynomial such that 
$$  \forall (a_1,\ldots, a_n)\in A: \ f(a_1,\ldots,a_n) =
p(a_1,\ldots,a_n).$$
\end{quote}
\end{Definition}

(There is a natural generalization of the ${\sigma}$-IP to $\kappa$-IP
 for any cardinal~$\kappa$. The theorems below can easily be
 generalized to these situations. 
   For the sake of simplicity we will not
 pursue this idea here.)

\begin{Theorem}
\label{mfl}
For every  lattice $L$ there is an infinite lattice $\bar L$
with the ${\sigma}$-IP, where  $L$ is a sublattice of~$\bar L$. 
\end{Theorem}

 We will sketch a proof of this theorem,
for simplicity
 only for the
case where the original lattice $L$ is of size $\le {\aleph_1}$. 

The simple (or even, by today's standards, trivial)
  set-theoretic idea behind this 
 theorem is the idea of transfinite iteration: 
\begin{quote}
{\sl If at first you don't succeed --- try, try again.}
\end{quote}

This idea goes back to Cantor.  The point here is, of course that
``again'' can mean here much more than ``many'' or even ``infinitely
many'' times.  When we are done with infinitely many repetitions, we
start over!

As an index set for this iteration we will use the well-ordered 
set\footnote{Some 
authors prefer to write $\Omega$ instead of~$\omega_1$.}
 $(\omega_1, <)$.   Rather than giving the set-theoretic definition
 of $\omega_1$ we will list some of its order-theoretic properties. 
The linear order  $ (\omega_1,{<})$
 can be characterized (up
  to isomorphism) by  (1)--(4):
  \begin{enumerate}
  \item $({\omega_1},{<})$ is a linear order.
  \item $({\omega_1},{<})$ is a well-order (i.e., every nonempty
  subset has a first element).
  \item ${\omega_1}$ is uncountable.
  \item For every $a\in {\omega_1}$ the set of predecessors,
  $\{i\in {\omega_1} : i < a\}$ is at most countable. (I.e., ``every bounded set
  is countable''.)
  \item Moreover: whenever $ A \subseteq {\omega_1} $ is countable,
  there there  is an $ a\in {\omega_1} $ such that  
  $   A \subseteq \{ i\in {\omega_1}: i  < a \}$. (I.e., ``every
  countable set is bounded.'')
  \end{enumerate}
  We  will write ${\aleph_1} $ for $|{\omega_1} |$.   We will write
$0$ for the smallest element of~$\omega_1$, and for $i\in
\omega_1 $ we will write $i+1$ for the ``successor'' of~$i$, i.e., 
$$ i+1 = \min \{j\in \omega_1: i < j \}.$$

It suffices 
to construct a sequence
 $(L_i:i\in
{\omega_1} )$ starting with the original lattice $L=L_0$  such that:  
\begin{enumerate}
\item[(A)] $ i < j $ implies: $L_i \subseteq L_j$ (as a sublattice).
\item[(B)] Each $L_i$ is countable. 
\item[(C)] For every $i \in {\omega_1}$, all natural numbers~$k$, 
and for all monotone partial functions with countable domain 
  $f:L_i^ k \to L_i$ there is a
 polynomial $p\in L_{i+1}[\x_1,\ldots, \x_k] $
  such that  for all $ \bar a = (a_1,\ldots,a_k) $ in the domain
  of~$f$: $f(\bar a)=p(\bar a)$.  
\end{enumerate}

We first explain why such a construction satisfying (A)+(B)+(C) 
is possible. 
 A rigorous argument would have to appeal to the fundamental
theorem of Definition by Transfinite Recursion; the following
informal argument captures the essence of the proof: 

Assume that this construction is not possible.  Since $\omega_1$
 is well-ordered this means that 
 there is a first index $j$ where the
construction breaks down.  If $j$ has no predecessor, then $L_j:=
\bigcup_{i<j} L_i$ satisfies all requirements, so we may assume that
$j$ has a predecessor, $j= i+1$. 

Given $L_i$, we now have to find $L_{i+1}$ such that every
countable partial monotone function from $L_i$ to $L_i$ is represented
by a polynomial with coefficients in~$L_{i+1}$.   Elementary cardinal
arithmetic shows that there are only ${\aleph_1} ^ {\aleph_0} =
{\aleph_1}  $  many such partial functions, so it is enough to apply
lemma~\ref{extend} repeatedly (${\aleph_1}$ many times)
to construct 
$L_{i+1}$. 

\medskip

We now explain why the properties (A)+(B)+(C) are sufficient: 

 $L:= \bigcup _{i\in {\omega_1} } L_i $ is a lattice, and (see
property (5) in the description of ${\omega_1} $) for every countable
	  partial function $ f:L \to L$ we can find an $i\in \omega_1$
 such that
both domain and range of $f$ are subsets of~$L_i$. Hence $f $ is
represented as a polynomial.

\section {Complete lattices}
\label{section.complete}

Recall that a lattice  is called ``complete''
 if every subset has a greatest lower and least
upper bound.  We say that a lattice is
$\sigma $-complete if every countable subset  has a greatest
lower and least upper bound. 
The lattice we constructed in theorem \ref{mfl} is easily seen to be
${\sigma} $-complete.   An important point in the construction of $\bar
L$ was the fact that every monotone partial function into a complete
lattice can be extended to a total monotone function, so completeness
(or $\sigma$-completeness, since the lattices we considered were
countable) played an important role in the proof of lemma 
\ref{extend}.  But how essential is completeness for
theorem~\ref{mfl} itself? 
Theorem \ref{sigma} below shows that there is really no
connection: $\sigma$-\opc\  is consistent with a strong negation of
$\sigma$-completeness.

\begin{Definition} Let $L$ be a lattice, $A, B \subseteq L$.  We write
$$ A < B$$ 
to abbreviate the property 
$$  \forall a\in A\, \forall b \in B:  a < b.$$
If $B$ is a singleton: $B = \{b\}$ then we may write $A < b$
instead of $A < \{b\}$, similarly if $A$ is a singleton. 
\end{Definition}

\begin{Definition} We call a bounded lattice $\sigma$-saturated if: 
\begin{quote}
Whenever $A, B \subseteq L$ are countable sets satisfying $A < B$,
$A$ upward directed\footnote{i.e., $\forall a,a'\in A \, \, 
\exists a^*\in A: a\le a^*, a'\le a^*$}, $B$ downward directed, 
\\
then there is an element $c\in L$ such that  $A <c < B$. 
\end{quote}
\end{Definition}

[This is a slight weakening of the usual model-theoretic notion of
saturation.]

In particular, in a $\sigma$-saturated lattice no countable increasing
sequence  has
a least upper bound. Thus, in all interesting cases (namely, in those
lattices in which there are infinite chains) 
 $\sigma$-saturation is a strong negation of
$\sigma$-completeness: In a ${\sigma} $-saturated lattice, no countable set
has  a least upper bound, except for trivial reasons.

 \begin{Theorem}\label{sigma} Let $L$ be a bounded lattice.  Then there is a
 $\sigma$-saturated $\sigma$-\opc
  lattice $\bar L$ which is an extension of~$L$. 
 \end{Theorem}

Consequently, $\sigma$-\opc\  does not imply $\sigma$-completeness. 

Using the method of transfinite iteration, it can be seen that 
theorem \ref{sigma} follows from lemma~\ref{extend}
together with the following lemma 
(which is a standard result in model theory, an easy consequence
of the compactness theorem):

\begin{Lemma} Let $L$ be a lattice, $A,B \subseteq L$, $A <
B$, $A$ upward directed, $B$ downward directed. 

Then there is a lattice $L^*$ extending $L$
 and an element $c\in L^*$
such that (in $L^*$) we have  $A <c < B$. 
\end{Lemma}

Note that even if we have $B = \{b\}$, $b = \sup A$ in $L$, $b$ will
lose this 
property in~$L^*$. Iterating this construction it is possible to
obtain a $\sigma$-saturated lattice~$\bar L$.

Interleaving this iteration with the construction from
 lemma~\ref{extend}
will guarantee that $\bar L$ will have the $\sigma$-IP.

\section{$\sigma$-polynomials}
\label{section.sigma}

In section \ref{section.opc}
 we have seen that there are no infinite lattices where
all monotone functions are represented by polynomials.  In section
\ref{section.IP}
 we have
seen that relaxing ``represented'' to ``${\sigma}$-interpolated'' we
get many lattices with this property, i.e., although there are no
infinite lattices which are \opc,  there are many lattices with the
${\sigma}$-IP. 

In this section we will consider another variation of this theme.  We
are now again interested in
{\em  representation} (rather than {\em interpolation}), 
 but instead of polynomial
functions we will consider a wider class of functions, by allowing
the infinitary operations $\sup
$ and~$\inf$.   We will take a ``conservative''   approach and only
consider $\inf$ and $\sup$ over countable sets. 

\begin{Definition}  Let $L$ be a lattice, $k\in \{1,2,\ldots\}$.   By 
$$ L[[\x_1,\ldots, \x_k]]$$
(the formal ``$\sigma $-polynomials'' over $L$) 
we denote the smallest set $S$  of formal expressions 
(in the formal variables 
$\x_1,\ldots, \x_k$) 
satisfying the
following: 
\begin{itemize}
\item The formal expressions $\x_1$,\dots, $\x_k$ are in~$S$, and all
elements of $L$ are in~$S$. 
\item Whenever $I$ is a finite or countably infinite set, and whenever
$(e_i:i \in I)$ is a family of expressions $e_i\in S$, 
then  also the formal
expressions 
$$ \sup (e_i: i \in I) \qquad \qquad
\inf (e_i: i \in I)
$$
are  in~$S$.
\end{itemize}
If $I= \{0,1\}$ then we may write $e_0\vee e_1$ and $e_0 \wedge e_1$
instead of 
$ \sup (e_i: i \in \{0,1\}) = \sup(e_0,e_1)$ and 
$\inf (e_i: i \in \{0,1\} ) = \inf(e_0,e_1)$, respectively. 

{\small 
[Several variants of this definition are possible --- we could divide
the set of formal expressions defined above by the set of all
equations that hold in all lattice extensions of $L$, or by the set
of equations that hold in all $\sigma$-complete lattice extensions of
$L$, etc.]}

Correspondingly, we let the set of $k$-ary  ``$\sigma$-polynomial
functions'' be the smallest set of (partial) functions from $L^k$ to $L$
containing all the projections and closed under the pointwise $\sup$
and $\inf $ over at most countable sets. 
\end{Definition}
There is a natural map (called ``pointwise evaluation'') from the set
of (formal)  $\sigma$-polynomials onto the set of  $\sigma$-polynomial
functions.

Clearly, all $\sigma$-polynomial functions are monotone. In general a
$\sigma$-polynomial function will not be total.

\begin{Theorem}
\label{s-p}
For any lattice $L$ there is a lattice $\bar L$ such that  $L$ is a
sublattice of $\bar L$, and $\bar L$ satisfies the following
properties: 
\begin{enumerate}
\item For every~$k$, every monotone function from $L^k$ to $L$ is a
$\sigma$-polynomial function, i.e., induced by a $\sigma$-polynomial.
\item All $\sigma$-polynomial functions on $\bar L$ are total.   Thus,
the monotone functions are exactly the ${\sigma}$-polynomial
functions. 
\item Moreover, $\bar L$ is complete, that is:  every
subset of $L$ has a greatest lower and a least upper bound. 
\end{enumerate}
\end{Theorem}
\begin{proof}
By  lemma \ref{extend} we can find a sequence $ L= L_0 \le L_1 \le \cdots $
of complete lattices such that:
\begin{quote}
For all~$k$, 
every monotone $f:L_n^ k\to L_n$ is represented by a polynomial with
coefficients in $L_{n+1} $. Moreover, $L_{n+1} $ is an end extension
of~$L_n$. 
\end{quote}
Let $L =\bigcup_n L_n$. It is easy to see
that $L$ is complete. {\small [If $A \subseteq L$, let $a_n:=
\inf_{L_n} (A\cap L_n)$.  Since the sequence $(a_n:n\in \N)$ is weakly
decreasing, and since $L_{n+1} $ is an end extension of~$L_n$, there
must be some $n_0$ such that   $\forall n\, a_n\in L_{n_0}$.  Now let
$a ^* := 
\inf\nolimits_{L_{n_o}}(a_n: n \in \N)
$, then clearly $a^* = \inf A$. ]}

Let $f:L^ k\to L$.  For simplicity assume~$k=1$. 
\\
For each $n$ define $f_n:L_n\to L_n$ by 
$$ f_n(x) := \sup\nolimits_{L_n} \{ y\in L_n: y\le_L f(x)\}.$$
$f_n$ is total (since $L_n$ is complete), and clearly monotone.  Let
$p_n\in L_{n+1}[\x]$ be a polynomial such that  $p_n(a) = f_n(a)$ for
all $a\in L_n$.  \\
Note: if $a,f(a) \in L_n$, then $p_n(a) = f(a)$. 

Now define a $\sigma $-polynomial $p(\x)$ by
$$ p(\x)= \mathop{\smash{\sup}\vrule width 0pt depth1pt}_{k\ge 0}\, 
          \mathop{\smash{\inf}\vrule width 0pt depth1pt}_{n\ge k} p_n(\x)$$
or slightly less formally: 
$$ p(\x)=
\bigl( p_0(\x) \wedge p_1(\x) \wedge p_2(\x) \wedge\cdots\bigr)
\ \vee \ 
\bigl(p_1(\x) \wedge p_2(\x) \wedge\cdots\bigr)
 \ \vee \ 
\bigl(p_2(\x) \wedge\cdots\bigr)
 \ \vee \ 
\cdots $$
It is easy to check that for all $a\in L$ the sequence
$(p_1(a), p_2(a), p_3(a), \ldots\, )$ is eventually constant with
value $f(a)$, and this implies $p(a) = f(a)$.
\end{proof}

Note that the lattices $L_n$ in this construction might have bigger
and bigger cardinalities.   Of course there are also many smaller
lattices in which every monotone function is a ${\sigma}$-polynomial,
e.g., every countable lattice with the IP. 
This suggests the following question: 

\begin{Problem}
\label{C}
  Describe all lattices in which every monotone 
function is a ${\sigma}$-polynomial function.  In particular, what are
 the cardinalities of such lattices? 
\end{Problem}

\section{Ortholattices}
\label{section.ortho}

This section contains a preview of results (without proofs)
 that will be published elsewhere (see \cite{ortho}).
  Again we give a theorem whose
 proof combines algebraic ideas with the simple set-theoretic idea of
 transfinite iteration.

An ortholattice is a bounded  lattice with an additional unary
operation $^\perp$ of ``orthocomplement'' which is antimonotone
 satisfying the following
laws: $x^ {\perp\perp} = x$, $x\vee x^ \perp = 1$, $x\wedge x^ \perp =
0$. (This implies de Morgan's laws. 
Informally, an ortholattice is a Boolean algebra without
distributivity.) 

Orthopolynomials and orthopolynomial functions 
are defined naturally:  in addition to the operations
$\vee$ and $\wedge$ we also allow $^ \perp$.   Hence orthopolynomials
are in general not monotone.    Is there any law that orthopolynomial
functions have to satisfy? The  theorem and its corollary 
 below  say that
orthopolynomials can have arbitrary behavior. 

We say that an ortholattice $\O$ has the ${\sigma}$-IP if every
(not necessarily monotone)
function $f:\O^ n\to \O$ is interpolated by an orthopolynomial on any
finite or countable set, and we say that an ortholattice $\O$ is
${\sigma}$-polynomially complete if every function $f:\O^ n\to \O$ is
represented by a ${\sigma}$-orthopolynomial.

\begin{Theorem}
\label{ortho}
For any ortholattice $(\O,\vee,\wedge,0,1,\perp) $ and any $f:\O\to\O$
there is an ortholattice $\bar\O$ such that  $\O$ is an
(ortho)sublattice of $\bar \O$, and there is an orthopolynomial
$p(\x) $ with coefficients in $\bar \O$ such that  $f(a) = p(a)$ for
all $a \in \O$. 
\end{Theorem}
\begin{Corollary}
 Every ortholattice $\O$ can be orthoembedded into an
ortholattice $\bar \O$ with the ${\sigma}$-IP. 
\end{Corollary}

\begin{proof}[Sketch]   

We only give a sketch of the main ideas of the proof.   The details
will appear elsewhere. 

Start with a lattice~$\O=\O_0$.  Let 
$$ \O_1 :=  \O_0 \, + \, (\O_0\times \O_0)$$
i.e., 
$ \O_1 $ is the ``horizontal sum'' of the disjoint  ortholattices 
$  \O_0 $ and $\O_0\times \O_0$, where we of course identify the two top
elements of the two lattices, and also the two bottom elements.

The functions $ g_1, g_2: \O \to \O\times \O $, defined by 
$$ \begin{array}{rcl}
g_1(x) &=& ( x, 0) \cr
g_2(x) &=& ( 1, x) \cr
   \end{array}
$$ 
are clearly  monotone.   Let $ A:= \{ (x,x^ \perp): x \in \O\}$.
This set is an antichain in~$\O_1$, so the functions 
$ h:A \to \O_0 \subseteq \O_1$, defined by 
$$ h(x,x^ \perp) = f(x)$$
is trivially a (partial) monotone function from $\O_1$ to~$\O_1$. 

Note that 
$$ f(x) = h(x,x^ \perp) = h ( \ g_1(x) \ \vee \ g_2(x)^ \perp \ ) $$
so $f$ can be written as a composition of monotone functions and the
orthocomplement function. 

By lemma \ref{extend}, there is a lattice extension $ L_2$ of the
lattice $ \O_1$ 
such that  $h,g_1,g_2$ can be represented by lattice polynomials with
coefficients in~$L_2$. 

A bit of work is now needed to find an ortholattice $ \O_3 $ which is
both a lattice extension of $L_2$ and an orthoextension of the
original ortholattice~$ \O$.  Once this ortholattice is found, the
function $ f $ can be represented as an orthopolynomial with
coefficients in $\O_3$.

\end{proof}

\end{document}